\documentclass[10pt]{article}

\usepackage{amssymb}
\usepackage{amsmath}

\newtheorem{theorem}{Theorem}

\begin{document}



\title
{Inverse boundary value problem for Schr\"odinger equation
in two dimensions}

\author{
M.~Yamamoto\thanks{ Department of Mathematical Sciences, University
of Tokyo, Komaba, Meguro, Tokyo 153, Japan e-mail:
myama@ms.u-tokyo.ac.jp}\, and \,O.~Yu.~Imanuvilov\thanks{ Department
of Mathematics, Colorado State University, 101 Weber Building, Fort
Collins, CO 80523-1874, U.S.A. E-mail: {\tt
oleg@math.colostate.edu}. This work is partially supported by NSF
Grant DMS 0808130.} }

\date{}

\maketitle
\begin{abstract}
We relax the regularity condition on potentials
of the Schr\"odinger equation in uniqueness results
on the inverse boundary value problem
which were recently proved in \cite{IUY} and \cite{Bu}.
\end{abstract}

Let $\Omega\subset \Bbb R^2$ be a bounded smooth domain with
$\partial\Omega=\cup_{j=0}^K\Sigma_j$ where $\Sigma_j$ are smooth
contours and $\Sigma_0$ is the external contour.  Let
$\nu=(\nu_1, \nu_2)$ be the unit outer normal to
$\partial\Omega$ and let $\frac{\partial}{\partial\nu} =
\nabla\cdot\nu$.

In this domain we
consider the Schr{\"o}dinger equation with some potential $q$:
\begin{equation}
(\Delta +q)u=0\quad\mbox{in}\,\,\Omega.
\end{equation}

Let $\widetilde \Gamma$ be a non-empty arbitrary fixed
relatively open subset of $\partial\Omega.$ Denote
$\Gamma_0=Int (\partial\Omega\setminus \widetilde \Gamma).$ Consider the
partial Cauchy data
\begin{equation}\label{popo}
\mathcal C_q= \left\{ \left(u,\frac{\partial u}{\partial \nu}\right)
\Biggl\vert_{\widetilde\Gamma}; \thinspace\thinspace
(\Delta+q)u=0\quad\mbox{in}\,\,\Omega,\,\,
u\vert_{\Gamma_0}=0, u\vert_{\widetilde\Gamma}=f\right\}.
\end{equation}

The goal of this article is  to improve the regularity assumption
on the potential $q$ in the case of arbitrary subboundary
$\widetilde{\Gamma}$ for the uniqueness result in the inverse
problem of recovery of potential from the partial data (\ref{popo}).
In the case of $\widetilde{\Gamma} = \partial\Omega$,
this inverse  problem  was formulated by  Calder\'on in
\cite{C}. Under the assumption $q\in C^{4+\alpha}(\overline\Omega)$
the result was proved in Imanuvilov, Uhlmann and Yamamoto
\cite{IUY}. In Guillarmou and Tzou \cite{GZ}, the assumption on
potentials was improved up to $C^{2+\alpha}(\overline\Omega).$

In particular, in the two-dimensional full Cauchy data case of
$\widetilde{\Gamma} = \partial\Omega$, we refer to Astala and
P\"aiv\"arinta \cite{AP},  Blasten \cite{EB},
Brown and Uhlmann \cite{B-U}, Bukhgeim \cite{Bu},
Nachman \cite{N}.  In \cite{EB}, the full Cauchy data uniquely
determine the potential within $W^1_p(\Omega)$ with $p>2$. As for
the related problem of recovery of the conductivity, \cite{AP}
proved the uniqueness result for conductivities from
$L^\infty(\Omega),$  improving the result of \cite{N}. We also
mention that for the case of full Cauchy data a relaxed regularity
assumption on potential  was claimed in \cite{Bu} but the proof
itself is missing some details.

In three or higher dimensions,
for the full Cauchy data, Sylvester and Uhlmann \cite{SU}
proved the uniqueness of recovery of conductivity in
$C^2(\overline\Omega)$, and later the regularity assumption was
relaxed up to $C^\frac 32(\overline\Omega)$ in P\"aiv\"arinta,
Panchenko and Uhlmann  \cite{PPU} and up to $W^\frac 32_p(\Omega)$
with $p>2n$ in Brown and Torres \cite{BT}.  For the case of partial
Cauchy data, uniqueness theorems were proved under assumption that a
potential of the Schr\"odinger equation belongs to
$L^\infty(\Omega)$ (see Bukhgeim and Uhlmann \cite{BuU}, Kenig,
Sj\"ostrand and Uhlmann \cite{KSU}).

Our main result is as follows
\begin{theorem}
Let  $q_1,q_2\in C^\alpha(\overline\Omega)$ for
some $\alpha\in (0,1)$ if $\widetilde \Gamma=\partial\Omega$
and $q_1,q_2\in W^1_p(\Omega)$ for some
$p>2$ otherwise. If $\mathcal C_{q_1}=\mathcal C_{q_2}$ then $q_1=q_2.$
\end{theorem}

The rest part of the paper is devoted to the proof of the
theorem.
Throughout the article, we use the following notations.

{\bf Notations.}
$i=\sqrt{-1}$, $x_1, x_2 \in
{\Bbb R}^1$, $z=x_1+ix_2$, $\overline{z}$
denotes the complex conjugate of $z \in \Bbb C$.
We identify $x = (x_1,x_2)
\in {\Bbb R}^2$ with $z = x_1 +ix_2 \in {\Bbb C}$.
$\partial_z = \frac 12(\partial_{x_1}-i\partial_{x_2})$,
$\partial_{\bar z}= \frac12(\partial_{x_1}+i\partial_{x_2})$,
$D = \left( \frac{1}{i}\partial_{x_1},
\frac{1}{i}\partial_{x_2}\right)$.
The tangential derivative on the boundary is given by
$\partial_{\vec\tau}=\nu_2\frac{\partial}{\partial x_1}
-\nu_1\frac{\partial}{\partial x_2}$, where
$\nu=(\nu_1, \nu_2)$ is the unit outer normal to
$\partial\Omega$.

{\bf Proof.}

{\bf First Step.}\\
Let $\Phi=\varphi+i\psi$ be a holomorphic function on
$\Omega$ such that $\varphi, \psi$ are real-valued and
\begin{equation}
\Phi\in C^2(\overline\Omega),\quad \mbox{Im}\, \Phi\vert_{\Gamma_0}=0.
\end{equation}
Denote by $\mathcal H$ the set of the critical
points of the function $\Phi.$ Suppose that this set is not empty,
each critical point is nondegenerate, $\mathcal
H\cap\overline\Gamma_0=\emptyset$ and
\begin{equation}\label{kk}
\mbox{mes}\, (\mathcal J)=0,\quad \mathcal
J=\{x; \thinspace
\partial_{\vec \tau}\psi(x)=0,x\in \widetilde \Gamma\}.
\end{equation}
Here $\vec \tau$ is an unit tangential vector to
$\partial\Omega.$   Consider the operator
$L_q(x,D)=-\sum_{j=1}^2(D_j+\tau i\varphi_{x_j})^2+q.$
It is known (see \cite{IUY1} Proposition 2.5)
that there exists a constant
$\tau_0$ such
that for $\vert \tau\vert\ge \tau_0$ and any $f\in L^2(\Omega)$,
there exists a solution to the boundary value problem
\begin{equation}\label{lola}
L_q(x,D)u =f\quad\mbox{in}\,\,\Omega, \quad u\vert_{\Gamma_0}=0
\end{equation}
 such that
 \begin{equation}\label{2}
 \Vert u\Vert_{H^{1,\tau}(\Omega)}/\root\of{\vert\tau\vert}\le C\Vert  f\Vert_{L^2(\Omega)}.
 \end{equation}
Moreover if $f/\partial_z\Phi\in L^2(\Omega)$,
then for any $\vert \tau\vert\ge \tau_0$  there exists
a solution to the boundary value problem (\ref{lola})
such that
\begin{equation}\label{3}
\Vert u\Vert_{H^{1,\tau}(\Omega)}
\le C\Vert  f/\partial_z\Phi\Vert_{L^2(\Omega)}.
\end{equation}
The constants $C$ in (\ref{2}) and (\ref{3}) are
independent of $\tau.$
Here and henceforth we set
$$
\Vert u\Vert_{H^{1,\tau}(\Omega)}
= (\Vert u\Vert_{H^1(\Omega)}^2
+ \vert \tau\vert^2\Vert u\Vert^2_{L^2(\Omega)})^{\frac{1}{2}}.
$$

{\bf Second Step.}\\
Here we will construct complex geometrical optics solutions.
Henceforth by $o_{L^2(\Omega)}(\frac 1\tau)$, we mean a function
$f(\epsilon, \tau, \cdot)\in L^2(\Omega)$ such that
$\lim_{\tau\to \infty}\vert\tau\vert
\Vert f(\epsilon,\tau,\cdot)\Vert_{L^2(\Omega)}
=0$ for all small $\epsilon > 0$, and
by $o(\frac 1\tau)$, we mean $a(\epsilon,\tau)$ such that
$\lim_{\tau\to \infty}\vert\tau\vert
\vert a(\epsilon,\tau)\vert =0$ for all small $\epsilon > 0$.

Let $\{q_{1,\epsilon}\}_{\epsilon\in(0,1)}$ be a sequence
of smooth functions converging to
$q_1$ in $W^1_p(\Omega)$ or $C^\alpha(\overline\Omega)$
(depending on the assumption on the regularity of $q_1$)
such that $q_{1,\epsilon}=q_1$ on $\mathcal H.$
Let $p_\epsilon $ be the complex geometrical optics solution to
the Schr\"odinger operator $\Delta+q_{1,\epsilon}$ which
we constructed in \cite{IUY}. The function $p_{\epsilon}$ can be
written in the form:
\begin{eqnarray}\label{mozilaal}
p_\epsilon(x)=e^{\tau{\Phi}}{(a+a_{0,\epsilon}/\tau)}
+e^{\tau\overline{\Phi}} \overline{(a+b_{1,\epsilon}/\tau)}\nonumber\\
- \bigg( e^{\tau\Phi}\frac{(\partial^{-1}_{\overline z}
(aq_{1,\epsilon})-M_{1,\epsilon})}{4\tau\partial_z\Phi} +
 e^{\tau\bar\Phi}\frac{(\partial^{-1}_{z}
(\overline{a}q_{1,\epsilon})-M_{3,\epsilon})}{4\tau\overline{\partial_z\Phi}}
\bigg )+e^{\tau\varphi}o_{L^2(\Omega)}(\frac 1\tau)\quad \mbox{as}\,\tau\rightarrow +\infty,
\end{eqnarray}
where $a\in C^6(\overline\Omega)$ is some holomorphic function
on $\Omega$ such that $\mbox{Re}\, a\vert_{\Gamma_0}=0$.
The operators
$\partial_{z}^{-1}$ and $\partial_{\overline z}^{-1}$
are given by
$$
\partial_{\overline z}^{-1}g=-\frac1\pi\int_\Omega
\frac{g(\zeta,\overline\zeta)}{\zeta-z}
d\xi_2d\xi_1,\quad
\partial_{ z}^{-1}g
= \overline{\partial^{-1}_{\overline{z}}\overline{g}},
$$
Moreover for some $\widetilde x\in \mathcal H$,
we assume that $a(\widetilde x) \ne 0$ and $a(x) = 0$ for
$x \in \mathcal H\setminus \{\widetilde x\}$,
and the polynomials
$M_{1,\epsilon}(z)$ and $M_{3,\epsilon}(\overline z)$ satisfy
$$
\partial_z^j(\partial^{-1}_{\overline z}(aq_{1,\epsilon})
-M_{1,\epsilon})(x)=0, \quad \quad
\partial^j_{\overline z}(\partial^{-1}_{z}(\overline{a}q_{1,\epsilon})
- M_{3,\epsilon})(x)= 0, \quad x \in \mathcal H,
$$
$a_{0,\epsilon},a_{1,\epsilon}\in C^6(\overline \Omega)$
are holomorphic functions such that
$$
(a_{0,\epsilon}+\overline{a}_{1,\epsilon})\vert_{\Gamma_0} =
\frac{(\partial^{-1}_{\overline z}
(aq_{1,\epsilon})-M_{1,\epsilon})}{4\partial_z\Phi}
+\frac{(\partial^{-1}_{z}
(\overline{a}q_{1,\epsilon})-M_{3,\epsilon})}
{4\overline{\partial_z\Phi}}.
$$

We look for a solution $u_1$ in the form $u_1=p_{\epsilon}+
m_{\epsilon}$. Consider the equation
$$
L_{q_1}(x,D)u_1=L_{q_{1, \epsilon}}(x,D)
(p_{\epsilon}+m_{\epsilon})
+(q_1- q_{1, \epsilon})(p_{\epsilon}+m_{\epsilon})
= L_{q_1}(x,D)m_{\epsilon}+(q_1- q_{1, \epsilon})p_{\epsilon}=0.
$$
By (\ref{3}) there exists a solution to the boundary value problem
$$L_{q_1}(x,D)m_{\epsilon}+(q_1- q_{1,
\epsilon})p_{\epsilon}=0\quad \mbox{in}\quad\Omega, \quad
m_{\epsilon}\vert_{\Gamma_0}=0
$$ such that
\begin{equation}\label{km}
\Vert m_{\epsilon}\Vert_{H^{1,\tau}(\Omega)}\le C(\epsilon)\quad
\forall \tau>\tau_0(\epsilon),
\end{equation}
where $C(\epsilon)$ is independent of $\tau$ and
$$
C(\epsilon)\rightarrow 0\quad \mbox{as}\quad
\epsilon\rightarrow 0.
$$

Since the Cauchy data (\ref{popo}) for potentials
$ q_1$ and $q_2$, are equal,
there exists a solution $u_2$  to the Schr\"odinger equation
with the potential $q_2$ such that $u_1=u_2$ on
$\partial\Omega$ and
$\frac{\partial u_1}{\partial\nu}
=\frac{\partial u_2}{\partial \nu}$
on $\widetilde \Gamma$.   Setting $u=u_1-u_2$, we obtain
\begin{equation}\label{pp}
(\Delta+q_2)u=(q_2-q_1)u_1\quad \mbox{in}\,\,\Omega, \quad
u\vert_{\partial\Omega}=\frac{\partial u}{\partial
\nu}\vert_{\widetilde\Gamma}=0.
\end{equation}

In a way similar to the construction of $u_1$,
we construct the complex geometrical optics solution
$v$ for the Schr\"odinger equation with the potential $q_2.$ The
construction of $v$ repeats the corresponding steps of the
construction of $u_1.$ The only difference is that instead of
$q_{1,\epsilon}$ and $\tau$, we use $q_{2,\epsilon}$ and $-\tau$
respectively. We provide details of the construction of $v$ for the
sake of completeness.

Let $\{q_{2, \epsilon}\}_{\epsilon\in (0,1)}$ be a sequence
of smooth functions converging to  sufficiently close to
$q_2$ in $W^1_p(\Omega)$ or $C^\alpha(\overline\Omega)$
such that $q_{2,\epsilon}=q_2$ on $\mathcal H.$
Let $\widetilde p_\epsilon $ be the complex geometrical
optics solution to the Schr\"odinger operator $\Delta+q_{2,\epsilon}$
constructed in \cite{IUY}:
\begin{eqnarray}\label{mozilaa}
\widetilde p_\epsilon(x)=e^{-\tau{\Phi}}{(a+b_{0,\epsilon}/\tau)}
+e^{-\tau\overline{\Phi}} \overline{(a+b_{1,\epsilon}/\tau)}\nonumber\\
+ \left (e^{-\tau\Phi}\frac{(\partial^{-1} _{\overline z}
(aq_{2,\epsilon})-M_{2,\epsilon})}{4\tau\partial_z\Phi}
+e^{-\tau\bar\Phi} \frac{(\partial^{-1}_{z}
(\overline{a}q_{2,\epsilon})-M_{4,\epsilon})}
{4\tau\overline{\partial_z\Phi}}\right
)+e^{-\tau\varphi}o_{L^2(\Omega)}(\frac 1\tau),
\end{eqnarray}
where $M_{2,\epsilon}(z)$ and $M_{4,\epsilon}(\overline z)$
satisfy
$$
\partial_z^j(\partial^{-1}_{\overline z}(aq_{1,\epsilon})
-M_{2,\epsilon})(x)=0, \quad \quad
\partial^j_{\overline z}(\partial^{-1}_{z}
(\overline{a}q_{1,\epsilon})
- M_{4,\epsilon})(x)= 0, \quad x \in \mathcal H.
$$
and $b_{0,\epsilon},b_{1,\epsilon}$ are holomorphic functions such
that
$$
(b_{0,\epsilon}+\overline
b_{1,\epsilon})\vert_{\Gamma_0}=-\frac{(\partial^{-1} _{\overline z}
(aq_{2,\epsilon})-M_{2,\epsilon})}{4\partial_z\Phi} -
\frac{(\partial^{-1}_{z}
(\overline{a}q_{2,\epsilon})-M_{4,\epsilon})}
{4\overline{\partial_z\Phi}}.
$$

We look for a solution $v$ in the form
$v=\widetilde p_{\epsilon}+ \widetilde m_{\epsilon}$. Consider the operator
$$
L_{q_2}(x,D)v=L_{q_{2, \epsilon}}(x,D)
(\widetilde p_{\epsilon}+\widetilde m_{\epsilon})
+(q_2- q_{2, \epsilon})(\widetilde p_{\epsilon}+
\widetilde m_{\epsilon})= L_{q_2}(x,D)\widetilde m_{\epsilon}
+(q_2- q_{2,\epsilon})\widetilde p_{\epsilon}=0.
$$
By (\ref{3}) there exists a
solution to the boundary value problem $$L_{q_2}(x,D)\widetilde
m_{\epsilon}+(q_2- q_{2, \epsilon})\widetilde p_{\epsilon}=0\quad
\mbox{in}\quad\Omega, \quad \widetilde m_{\epsilon}\vert_{\Gamma_0}=0
$$ such that
\begin{equation}\label{km1}
\Vert \widetilde m_{\epsilon}\Vert_{H^{1,\tau}(\Omega)}\le
C(\epsilon)\quad \forall \tau>\tau_0(\epsilon),
\end{equation}
where $C(\epsilon)$ is independent of $\tau$ and
$$
C(\epsilon)\rightarrow 0\quad \mbox{as}\quad \epsilon\rightarrow 0.
$$

{\bf Third Step.}
\\
We will prove $q_1(\widetilde{x})=q_2(\widetilde{x})$
where $a(\widetilde{x})\ne 0$ and $a(x) = 0$ for
$x \in \mathcal{H} \setminus
\{ \widetilde{x}\}$ in the case where $q_1, q_2 \in
W^1_p(\Omega)$.

Denote $q=q_1-q_2.$ Taking the scalar product of equation (\ref{pp})
and the function $v$, we have:
\begin{equation}
\int_\Omega q u_1vdx=0. \end{equation} By (\ref{km}) and (\ref{km1})
\begin{equation}\label{01}
0=\int_\Omega qu_1vdx=\int_\Omega q p_{ \epsilon} \widetilde
p_{\epsilon}dx+K( \epsilon,\tau),
\end{equation} where
\begin{equation}\label{op}
\overline{\lim_{\tau\rightarrow +\infty}}\tau\vert K(
\epsilon,\tau)\vert\le C(\epsilon), \quad C( \epsilon)\rightarrow
0\quad \mbox{as}\quad\epsilon \rightarrow 0.
\end{equation}
From (\ref{01}), (\ref{op}) and the explicit formulae
(\ref{mozilaal}), (\ref{mozilaa}) for the
construction of complex geometrical optics solutions, we have
$$
\int_\Omega q(a^2+\overline a^2)dx=0.
$$
Computing the remaining terms, we have:
\begin{eqnarray}\label{lala}
K(\epsilon,\tau)+ \frac 1\tau\int_\Omega q(
a(a_{0,\epsilon}+b_{0,\epsilon})+\overline{a(a_{1,\epsilon}+b_{1,\epsilon})})dx
+\int_\Omega
q(a\overline a e^{2\tau i\psi}+ a\overline a e^{-2\tau i\psi})dx \nonumber\\
+\frac{1}{4\tau}\int_{\Omega} \left( qa \frac{\partial_{\overline
z}^{-1}(a q_{2,\epsilon})-M_{2,\epsilon}} {\partial_z\Phi} +
q\overline{a}\frac{\partial_{z}^{-1}(q_{2,\epsilon}\overline{a})
-{M_{4,\epsilon}}}{\overline{\partial_z\Phi}}\right)dx\nonumber\\
-\frac{1}{4\tau}\int_\Omega\left( qa\frac{\partial_{\overline
z}^{-1}(q_{1,\epsilon}a)-M_{1,\epsilon}}{\partial_z\Phi} +q\overline
a\frac{\partial_{z}^{-1}(q_{1,\epsilon}\overline a)-{
M_{3,\epsilon}}}{\overline{\partial_z\Phi}}\right)dx\nonumber\\+
o(\frac{1}{\tau})=0\quad\mbox{as}\,\,\tau \rightarrow +\infty.
\end{eqnarray}

Since the functions $q_j$ are not supposed to be from
$C^2(\overline\Omega)$, we can not directly use the
stationary phase argument (e.g., Evans \cite{E}).
Consider two cases. Assume that $q\in
W^1_p(\Omega)$ with $p>2.$ We have
\begin{equation}\label{rono}
\int_\Omega q\mbox{Re}\,(a\overline a e^{2\tau i\psi})dx =\int_\Omega
q_{\epsilon}\mbox{Re}\,(a\overline a e^{2\tau i\psi})dx+\int_\Omega
(q-q_{\epsilon})\mbox{Re}\,(a\overline a e^{2\tau i\psi})dx.
\end{equation}
We set $q_{\epsilon}=q_{1,\epsilon}-q_{2,\epsilon}$.  Taking into
account that $q_{j,\epsilon}=q_j$ on $\mathcal H$, $j=1,2$, (\ref{kk})
and using the stationary phase argument, similar to \cite{IUY},
we compute
\begin{equation}\label{masa}
\int_\Omega q_{\epsilon}(a\overline a e^{2\tau i\psi}
+ a\overline a e^{-2\tau
i\psi})dx=\frac{2\pi (q\vert a\vert^2)(\widetilde
x)\mbox{Re}\,e^{2{\tau} i\mbox{Im}\,\Phi(\widetilde x)}} {{\tau}
\vert(\mbox{det}\thinspace \mbox{Im}\,\Phi'')(\widetilde
x)\vert^\frac
12}+o\left(\frac{1}{{\tau}}\right)\quad\mbox{as}\,\,\tau \rightarrow
+\infty.
\end{equation}
For the second integral in (\ref{rono}) we obtain
$$
\int_\Omega (q-q_{\epsilon})(a\overline a e^{2\tau i\psi}+ a\overline a
e^{-2\tau i\psi})dx=\int_\Omega (q-q_{ \epsilon})\left (a\overline a
\frac {(\nabla\psi,\nabla)e^{2\tau i\psi}}{2\tau i\vert
\nabla\psi\vert^2}-a\overline a  \frac {(\nabla\psi,\nabla)e^{-2\tau
i\psi}}{2\tau i\vert \nabla\psi\vert^2}\right )dx
$$
$$
= \int_{\partial\Omega} (q-q_{ \epsilon})\left (a\overline a \frac
{(\nabla\psi,\nu)e^{2\tau i\psi}}{2\tau i\vert
\nabla\psi\vert^2}-a\overline a  \frac {(\nabla\psi,\nu)e^{-2\tau
i\psi}}{2\tau i\vert \nabla\psi\vert^2}\right )d\sigma
$$
\begin{equation}\label{opl}-\frac{1}{2\tau i}
\int_\Omega\left \{ e^{2\tau
i\psi}\mbox{div}\,\left ((q-q_{\epsilon})a\overline a \frac
{\nabla\psi}{\vert \nabla\psi\vert^2}\right )- e^{-2\tau
i\psi}\mbox{div}\,\left ( (q-q_{ \epsilon})a\overline a \frac
{\nabla\psi}{\vert \nabla\psi\vert^2}\right )\right\}dx.
\end{equation}
Since $\psi\vert_{\Gamma_0}=0$ we have
$$
\int_{\partial\Omega} (q-q_{ \epsilon})a\overline a \left (\frac
{(\nabla\psi,\nu)e^{2\tau i\psi}}{2\tau i\vert \nabla\psi\vert^2}-
\frac {(\nabla\psi,\nu)e^{-2\tau i\psi}}{2\tau i\vert
\nabla\psi\vert^2}\right )d\sigma=\int_{\widetilde\Gamma}
\frac{(q-q_{
\epsilon})a\overline a}{2\tau i\vert \nabla\psi\vert^2}
(\nabla\psi,\nu)(e^{2\tau i\psi}- e^{-2\tau i\psi})d\sigma.
$$
By (\ref{kk}) and Proposition 2.4 in \cite{IUY} we have that
\begin{equation}\nonumber
\int_{\partial\Omega} (q-q_{ \epsilon})a\overline a \left (\frac
{(\nabla\psi,\nu)e^{2\tau i\psi}}{2\tau i\vert \nabla\psi\vert^2}-
\frac {(\nabla\psi,\nu)e^{-2\tau i\psi}}{2\tau i\vert
\nabla\psi\vert^2}\right )d\sigma=o(\frac 1\tau)\quad
\mbox{as}\,\,\tau\rightarrow+\infty.
\end{equation}
The last integral over $\Omega$ in
formula (\ref{opl}) is $o(\frac{1}{\tau})$ and so
\begin{equation}\label{-3}
\int_\Omega(q-q_{ \epsilon})(a\overline a e^{2\tau i\psi}
+ a\overline a e^{-2\tau
i\psi})dx=o(\frac 1\tau)\quad\mbox{as}\,\,\tau \rightarrow +\infty.
\end{equation}
Taking into account that $\psi(\widetilde x)\ne 0$ and using (\ref{masa}), (\ref{-3}) we have
from (\ref{lala}) that
\begin{equation}
\frac{2\pi (q\vert a\vert^2)(\widetilde x)} {
\vert(\mbox{det}\thinspace \mbox{Im}\,\Phi'')(\widetilde
x)\vert^\frac 12}+\widetilde C(\epsilon)=0,
\end{equation}
where $\widetilde C(\epsilon)\rightarrow +0$ as $\epsilon\rightarrow 0.$
Hence
\begin{equation}\label{victory}
q(\widetilde x)=0 \quad\mbox{if $a(\widetilde{x})\ne 0$ and
$a(x) = 0$ for $x \in \mathcal{H}\setminus \{ \widetilde x\}$}.
\end{equation}
Since a point $\tilde x$ can be chosen arbitrarily close to any
given point in $\Omega$ (see \cite{IUY}), we have $q\equiv 0$,
that is, the proof of the theorem is completed if
$q_1, q_2 \in W^1_p(\Omega)$.

{\bf Fourth Step.}
\\
Now let $q\in C^\alpha(\overline\Omega)$ with
some $\alpha\in (0,1)$ and
$\partial\Omega=\widetilde \Gamma.$

We recall the following classical result of H\"ormander \cite{Ho}.
 Consider the
"oscillatory integral operator"
$$
T_\tau f(x)=\int_\Omega e^{-\tau i\psi(x,y)}a(x,y)f(y)dy,
$$
where $\psi\in C^\infty(\Bbb R^2\times \Bbb R^2)$ and
$a(\cdot,\cdot)\in C_0^\infty(\Bbb R^2\times \Bbb R^2).$
We introduce the following matrix
$$
H_\psi=\{\partial^2_{x_iy_j}\psi\}.
$$
\begin{theorem}\label{opop} Suppose that $det\, H_\psi\ne 0$ on
$\mbox{supp}\, a.$
Then
$$
\Vert T_\tau\Vert_{L^2\rightarrow L^2}\le \frac{C}{\tau}.
$$
\end{theorem}

Consider our holomorphic function
$\Phi(x,y)=(x_1+ix_2-(y_1+iy_2))^2+i.$ We set
$\psi(x,y)=2(x_1-y_1)(x_2-y_2)-1.$  Then
$$
H_\psi(x,y)=\left(
\begin{array}{cccc}
0 & -2\\
-2 & 0\\
\end{array}\right)
$$
and $det H_\psi(x,y)=-4.$  Then the condition in Theorem \ref{opop}
holds true.

We set $a(x,y)=\chi(x)\chi(y)$ where $\chi\in C_0^\infty(\Bbb R^n)$
and $\chi\vert_\Omega\equiv 1.$ Then, by  Theorem \ref{opop},  there
exists a constant $C$ independent of $\tau$ such that
\begin{equation}\label{-2}
\Vert T_\tau\Vert_{L^2\rightarrow L^2}+\Vert
T_{-\tau}\Vert_{L^2\rightarrow L^2}\le C/\tau.
\end{equation}
Setting $f=(q-q_{ \epsilon})a\overline a\chi_\Omega$ by (\ref{-2}) we have
\begin{equation}\label{-1}
\Vert T_\tau f\Vert_{L^2(\Omega)}+\Vert
T_{-\tau}f\Vert_{L^2(\Omega)}\le C(\epsilon)/\tau, \quad
C(\epsilon)\rightarrow 0\quad\mbox{as}\quad \epsilon
\rightarrow +0.
\end{equation}
Therefore, by (\ref{-1}), in the ball $B(\widetilde x,\delta) \equiv
\{ x;\thinspace \vert x - \widetilde{x}\vert < \delta\}$,
there exists a sequence
of points $y(\tau)$  such that
\begin{equation}\label{moso}
\vert (T_\tau) f(y(\tau))\vert
+ \vert (T_{-\tau})f(y(\tau))\vert
\le \frac{C\epsilon}{\tau\delta^2}.
\end{equation}
Let $y(\tau)=(y_1(\tau),y_2(\tau))\rightarrow \hat y(\epsilon)$ as
$\tau\rightarrow +\infty.$ By the stationary phase argument taking
into account that $\psi(\widetilde x,\widetilde x)=-1$, we have
\begin{equation}\label{masa}
\int_\Omega (q_\epsilon-(q_\epsilon-q)(y(\tau))\mbox{Re}
\{a\overline a e^{-2\tau i \psi(y(\tau),x)}\}dx=\frac{2\pi (q\vert
a\vert^2)(\hat y(\epsilon))\mbox{Re}\,e^{2{\tau}i }} {{\tau}
}+o\left(\frac{1}{{\tau}}\right).
\end{equation}
From (\ref{lala}), (\ref{masa}), (\ref{moso}) we obtain
\begin{equation}
{2\pi (q\vert a\vert^2)(\hat y(\epsilon)) \mbox{Re}\,e^{2{\tau} i }}
+\widetilde C(\epsilon)=0,
\end{equation}
where $\overline{\lim}_{\tau\rightarrow+\infty}\vert\widetilde
C(\epsilon)\vert\rightarrow +0$ as $\epsilon\rightarrow 0.$
Therefore as $\epsilon$ goes to zero, we have
$$
q(\hat x)=0.
$$
Here $\hat x \in B(\widetilde x,\delta)$ such
that $\hat y(\epsilon)\rightarrow \hat x$ as
$\epsilon \rightarrow +0$.
Since $\delta>0$ and $\widetilde x$ are chosen arbitrarily,
we conclude that $q\equiv 0$ in $\Omega$.
Thus the proof of the theorem is completed.
$\square$

\end{document}